\documentclass[reqno,11pt]{amsart}

\usepackage{amsthm,amssymb,amsmath,amscd,graphicx,latexsym,enumerate}

\newcommand{\CC}{{\mathbb C}}

\newcommand{\RR}{{\mathbb R}}
\renewcommand{\SS}{{\mathbb S}}

\newcommand{\K}{{\mathcal K}}

\renewcommand{\O}{{\mathcal O}}

\newcommand{\fg}{{\mathfrak g}}
\newcommand{\fh}{{\mathfrak h}}
\newcommand{\ft}{{\mathfrak t}}
\newcommand{\fk}{{\mathfrak k}}

\newcommand{\fp}{{\mathfrak p}}

\newtheorem{theorem}{Theorem}[section]

\newtheorem{definition}{Definition}[section]

\newtheorem{proposition}{Proposition}[section]

\theoremstyle{definition}
\newtheorem*{theorem*}{Theorem}
\newtheorem*{proposition*}{Proposition}
\newtheorem*{definition*}{Definition}
\newtheorem*{remark*}{Remark}
\newtheorem*{remarks*}{Remarks}
\newtheorem*{lemma*}{Lemma}
\newtheorem{example*}{Example}

\input xy
\xyoption{all}

\numberwithin{equation}{section}

\title[Homogeneous Spaces]{A GKM description of the equivariant cohomology ring of a homogeneous space}
\author[V. Guillemin]{V. Guillemin}
\address{Department of Mathematics, MIT, Cambridge, MA 02139}
\email{vwg@math.mit.edu}
\author[T. Holm]{T. Holm}
\address{Department of Mathematics, MIT, Cambridge, MA 02139}
\email{tsh@math.mit.edu}
\author[C. Zara]{C. Zara}
\address{Department of Mathematics, Yale University, New Haven, CT 06520}
\email{zara@math.yale.edu}
\date{\today}

\begin{document}

\begin{abstract}
Let $T$ be a torus of dimension $n>1$ and $M$ a compact
$T-$manifold. $M$ is a GKM manifold if the set of zero dimensional
orbits in the orbit space $M/T$ is zero dimensional and the set of
one dimensional orbits in $M/T$ is one dimensional. For such a
manifold these sets of orbits have the structure of a labelled
graph and it is known that a lot of topological information about $M$
is encoded in this graph.

In this paper we prove that every compact homogeneous space $M$ of
non-zero Euler characteristic is of GKM type and show that the
graph associated with $M$ encodes \emph{geometric} information
about $M$ as well as topological information. For example, from
this graph one can detect whether $M$ admits an invariant complex
structure or an invariant almost complex structure.
\end{abstract}

\maketitle

\tableofcontents

\section{Introduction}
\label{sec:1}

Let $T$ be a torus of dimension $n>1$, $M$ a compact manifold,
$$\tau : T \times M \to M$$
a faithful action of $T$ on $M$, and $M/T$ the orbit space of
$\tau$. $M$ is called a \emph{GKM manifold} if the set of zero
dimensional orbits in the orbit space $M/T$ is zero dimensional
and the set of one dimensional orbits in $M/T$ is one dimensional.
Under these hypotheses, the union, $\Gamma \subset M/T$, of the
set of zero and one dimensional orbits has the structure of a
graph: Each connected component of the set of one-dimensional
orbits has at most two zero-dimensional orbits in its closure; so
these components can be taken to be the edges of a graph and the
zero-dimensional orbits to be the vertices. Moreover, each edge,
$e$, of $\Gamma$ consists of orbits of the same orbitype: namely,
orbits of the form ${\mathcal O}_e=T/{H_e}$, where $H_e$ is a
codimension one subgroup of $T$. Hence one has a labelling
\begin{equation} \label{eq:1.1}
e \to H_e
\end{equation}
of the edges of $\Gamma$ by codimension one subgroups of $T$.

It has recently been discovered that if $M$ has either a
$T-$invariant complex structure or a $T-$invariant symplectic
structure, the data above - the graph $\Gamma$ and the labelling
\eqref{eq:1.1} - contain a surprisingly large amount of
information about the global topology of $M$. For instance,
Goresky, Kottwitz, and MacPherson proved that the ring structure
of the equivariant cohomology ring $H_T^*(M)$ is completely
determined by this data, and Knutson and Rosu have shown that the
same is true for the ring $K_T(M)\otimes \CC$.

The manifolds $M$ which we will be considering below will be
neither complex nor symplectic; however we will make an assumption
about them which is in some sense much stronger then either of
these assumptions. We will assume that $T$ is the Cartan subgroup
of a compact, semisimple, connected Lie group $G$, and that $G$
acts transitively on $M$, \emph{i.e.} $M$ is a $G$-homogeneous
space. There is a simple criterion for such a manifold to be a GKM
manifold.

\begin{theorem}\label{th:new1.1}
Suppose $M$ is a $G$-homogeneous manifold.  Then the following are
equivalent.
\begin{enumerate}
\item The action of $T$ on $M$ is a GKM action;
\item The Euler characteristic of $M$ is non-zero;
\item $M$ is of the form $M=G/K$, where $K$ is a closed
subgroup of $G$ containing $T$.
\end{enumerate}
\end{theorem}

As we mentioned above, the data \eqref{eq:1.1} determine the ring
structure of $H^*_T(M)$ if $M$ is either complex or symplectic. This
result is, in fact, true modulo an assumption which is weaker than
either of these assumptions; and this assumption - equivariant
formality - is satisfied by homogeneous spaces which satisfy the
hypotheses of the theorem. Hence, for these spaces, one has two
completely different descriptions of the ring $H_T^*(M)$: the
graph theoretical description above and the classical Borel
description. In Section~\ref{sec:2} we will compute the graph
$\Gamma$ of a space $M$ of the form $G/K$, with $T \subset K$, and
show that it is a \emph{homogeneous} graph, i.e. we will show that
the Weyl group of $G$, $W_G$, acts transitively on the vertices of
$\Gamma$ and that this action preserves the labelling
\eqref{eq:1.1}. We will then use this result to compare the two
descriptions of $H_T^*(M)$.

One of the main goals in this paper is to show that for
homogeneous manifolds $M$ of GKM type, some important features of
the geometry of $M$ can be discerned from the graph $\Gamma$ and
the labelling \eqref{eq:1.1}. One such feature is the existence of
a $G-$invariant almost complex structure. The subgroups, $H_e$,
labelling the edges of $\Gamma$ are of codimension one in $T$; so,
up to sign, they correspond to weights, $\alpha_e$, of the group
$T$. It is known that the $W_K-$invariant labelling \eqref{eq:1.1}
can be lifted to a $W_K-$invariant labelling
\begin{equation}\label{eq:1.2}
e \to \alpha_e
\end{equation}
if $M$ is a coadjoint orbit of $G$ (hence, in particular, a
complex $G-$manifold). Moreover, this labelling has certain simple
properties which we axiomatize by calling a map with these
properties an \emph{axial function} (see Section~\ref{ssec:3.1}).
In Section~\ref{sec:3} we prove the following result.

\begin{theorem}\label{th:new1.2}
The homogeneous space $M$ admits a $G-$invariant almost complex
structure if and only if $\Gamma$ possesses a $W_K-$invariant
axial function \eqref{eq:1.2} compatible with \eqref{eq:1.1}.
\end{theorem}

This raises the issue: Is it possible to detect from the graph
theoretic properties of the axial function \eqref{eq:1.2} whether
or not $M$ admits a $G-$invariant complex structure? Fix a vector
$\xi \in \ft$ such that $\alpha_e(\xi) \neq 0$ for all oriented
edges, $e$, of $\Gamma$, and orient these edges by requiring that
$\alpha_e(\xi)>0$. We prove in Section \ref{sec:4} the following
theorem.

\begin{theorem}\label{th:new1.3}
A necessary and sufficient condition for $M$ to admit a
$G-$invariant complex structure is that there exist no oriented
cycles in $\Gamma$.
\end{theorem}

\noindent {\bf Remarks:}

\begin{enumerate}
\item $M$ admits a $G-$invariant complex structure if and only if
it admits a $G-$invariant symplectic structure; and, by the
Konstant-Kirillov theorem, it has either (and hence both) of these
properties if and only if it is a coadjoint orbit of $G$.

\item By the Goresky-Kottwitz-MacPherson theorem, the graph $\Gamma$
and the axial function \eqref{eq:1.2} determine the cohomology
ring structure of $M$. The \emph{additive} cohomology of $M$, i.e.
its Betti numbers, $\beta_i$, can be computed by the following
simple recipe: For each vertex, $p$, of the graph $\Gamma$, let
$\sigma_p$ be the number of oriented edges issuing from $p$ with
the property that $\alpha_e(\xi)<0$. Then
$$\beta_i =
  \begin{cases}
    0, & \text{ if $i$ is odd}, \\
    \#\{p; \sigma_p=i/2\}, & \text{ if $i$ is even}.
\end{cases} \;
$$
\end{enumerate}

One question we have not addressed in this paper is the question:
When is a labelled graph the GKM graph of a homogeneous space of
the form $G/K$ with $T \subset K$? For some partial answers to this
question see \cite{Ho}.

{\bf Acknowledgements:} We are grateful to David Vogan for helping
us formulate and prove the results in Section \ref{sec:3} and to
Bert Kostant for pointing out to us that a homogeneous space is a
quotient of a compact group by a closed subgroup of the same rank
if and only if its Euler characteristic is non-zero, and for
making us aware of a number of nice properties of such spaces.

\section{The Equivariant Cohomology of Homogeneous Spaces}
\label{sec:2}

\subsection{The Borel description}
\label{ssec:2.1}

Let $G$ be a compact semi-simple Lie group, $T$ a Cartan
subgroup of $G$, $K$ a closed subgroup of $G$ such that
$$T \subset K \subset G \; , $$
and let $\ft \subset \fk \subset \fg$ be the Lie algebras of $T$,
$K$, and $G$.

Let $\Delta_K \subset \Delta_G$ be the roots of $K$ and $G$, with
$\Delta_K^+ \subset \Delta_G^+$ sets of positive roots, let
$$\Delta_{G,K} = \Delta_G  - \Delta_K \; , $$
and let $W_K \subset W_G$ be the Weyl groups of $K$ and $G$. We
will regard an element of $W_G$ both as an element of $N(T)/T$ and as
a transformation of the dual Lie
algebra $\ft^*$ (or as a transformation of $\ft$, via the isomorphism
$\ft^* \simeq \ft$ given by the Killing form).
Also, we will assume for simplicity that $G$ is simply connected and
that the homogeneous
space $G/K$ is oriented.


Now suppose $M$ is a $G$-manifold. Then the  equivariant cohomology
ring $H^*_T(M)$ is related to the cohomology ring $H^*_G(M)$ by
$$
  H^*_T(M) = H^*_G(M) \otimes_{\SS(\ft^*)^{W_G}} \SS(\ft^*) \; .
$$
(see \cite[Chap. 6]{GS}), where $\SS(\ft^*)$ is the symmetric
algebra of $\ft^*$. In particular, let $M=G/K$, where $K$ acts on
$G$ by \emph{right} multiplication. Then $G$ acts on $M$ by
\emph{left} multiplication and
$$
  H^*_G(M) = H^*_G(G/K) = \SS(\fk^*)^K = \SS(\ft^*)^{W_K} \; ,
$$
hence
\begin{equation}\label{eq:2.1}
  H^*_T(G/K) = \SS(\ft^*)^{W_K} \otimes_{\SS(\ft^*)^{W_G}}
  \SS(\ft^*)\; .
\end{equation}
This is the Borel description of $H^*_T(G/K)$. Throughout
this paper, unless stated otherwise, $M$ is the homogeneous space
$G/K$.

\subsection{The GKM graph of $M$}
\label{ssec:2.2}

In the following subsections, we will show that $M$ is equivariantly
formal and is a GKM space.  Then we will relate the the GKM
description of the equivariant cohomology ring of $M$ to the
description above.

\subsubsection{Equivariant formality}

The $\SS(\ft^*)-$module structure of the equivariant cohomology
ring $H^*_T(M)$ can be computed by a spectral sequence (see
\cite[p. 70]{GS}) whose $E_2$ term is $H(M)\otimes \SS(\ft^*)$,
and if this spectral sequence collapses {\it at this stage}, then
$M$ is said to be \emph{equivariantly formal}. If $M=G/K$, with $T
\subset K$, then,
$$H^{odd}(M) = 0, $$
(see \cite[p. 467]{GHV}), and from this it is easy to see that all
the higher order coboundary operators in this spectral sequence
have to vanish by simple degree considerations. Hence $M$ is
equivariantly formal. One implication of equivariant formality is
the following:

\begin{theorem}
The restriction map
\begin{equation}\label{eq:new1}
H^*_T(M) \to H^*_T(M^T)
\end{equation}
is injective.
\end{theorem}

\begin{proof}
By a localization theorem of Borel (see \cite{Bo}), the kernel of
\eqref{eq:new1} is the torsion submodule of $H_T^*(M)$. However,
if $M$ is equivariantly formal, then $H^*_T(M)$ is free as an
$\SS(\ft^*)-$module, so the kernel has to be zero.
\end{proof}

Thus $H^*_T(M)$ imbeds as a subring of the ring
\begin{equation}\label{eq:new2}
H^*_T(M^T)=H^0(M^T) \otimes_{\CC} \SS(\ft^*) \; .
\end{equation}

\noindent We will give an explicit description of this subring in
Section~\ref{ssec:2.3}.

\subsubsection{The Euler characteristic.} It follows from
\eqref{eq:new2} that, if $M$ is a homogeneous space of the form
$G/K$, with $T \subseteq K$, then the Euler characteristic of $M$
is equal to
$$\chi(M) = \sum_{i} \dim H^{2i}(M) \; ;$$
in particular, the Euler characteristic is non-zero. It is easy to
see that the converse is true as well.

\begin{proposition}
If $M=G/K$ and the rank of $K$ is strictly less than the rank of
$G$, then the Euler characteristic of $G/K$ is zero.
\end{proposition}

\begin{proof}
Let $h$ be an element of $T$ with the property that
$$\{ h^N \; ; -\infty < N < \infty \}$$
is dense in $T$. Suppose that the action of $h$ on $G/K$ fixes a
coset $g_0K$. Then $g_0^{-1}hg_0 \in K$, i.e. $h$ is conjugate to
an element of $K$ and hence conjugate to an element $h_1$ of the
Cartan subgroup $T_1$ of $K$. However, if the iterates of $h$ are
dense in $T$, so must be the iterates of $h_1$ and hence $T_1 =
T$. Suppose now that $h=\exp{\xi}, \xi \in \ft$. If $h$ has no
fixed points, then the vector field $\xi_M$ can have no zeroes and
hence the Euler characteristic of $M$ has to be zero.
\end{proof}

\subsubsection{The fixed points}

We prove in this section that the action of $T$ on $M$ is a GKM
action; i.e. that the set of zero dimensional orbits in the orbit
space $M/T$ is zero dimensional, and the set of one dimensional
orbits is one dimensional. It is easy to see that these properties
are equivalent to

\begin{itemize}
\item[(1)] $M^T$ is finite;
\item[(2)] For every codimension one subgroup $H$ of $T$, $\dim M^H
\leq 2$.
\end{itemize}

We will show that if $M$ is of the form $G/K$, with $T \subseteq K$,
then it has the two properties above, and we will also show that
it has the following third property:

\begin{itemize}
\item[(3)] For every subtorus $H$ of $T$ and every connected
component $X$ of $M^H$, $X^T \neq \emptyset$.
\end{itemize}

It is well known that these properties hold for the homogeneous
space $\O = G/T$. The first two properties can be checked directly
(see \cite{GZ1}, and the third property holds because $\O$ is a
compact symplectic manifold and the action of $T$ is Hamiltonian.
Therefore, to prove that $M$ satisfies properties 1-3, it suffices
to prove the following theorem.

\begin{theorem}\label{th:new2}
For every subtorus $H$ of $T$, the map
\begin{equation}\label{eq:new3}
\O =G/T \to G/K = M
\end{equation}
sends $\O^H$ onto $M^H$.
\end{theorem}

\begin{proof}
Let $p_0$ be the identity coset in $M$ and $q_0$ the identity
coset in $\O$. Let $h$ be an element of $H$ with the property that
$$\{ h ^N \; ; \; -\infty < N < \infty \}$$
is dense in $H$. If $p=gp_0 \in M^H$, then $g^{-1}hg \in K$; so
$g^{-1}hg = ata^{-1}$, with $a \in K$ and $t \in T$. Thus
$hga=gat$ and hence $hq=q$, where $q=gaq_0$. But under the map
\eqref{eq:new3}, $q_0$ is sent to $p_0$, so $q$ is sent to
$gap_0=gp_0=p$. \end{proof}

In particular, Theorem~\ref{th:new2} tells us that the map $\O^T
\to M^T$ is surjective. However,
$$\O^T = N_G(T)/T = W_G \; , $$
so $M^T$ is the image of $W_G = N_G(T)/T$ in $G/K$. But $N_G(T)
\cap K = N_K(T)$, the normalizer of $T$ in $K$, so
$$(N_G(T) \cap K)/T = W_K \; , $$
and hence we proved:

\begin{proposition}
There is a bijection
$$M^T \simeq W_G/W_K \; ;$$
in particular, $W_G =N_G(T)/T$ acts transitively on $M^T$.
\end{proposition}

\subsubsection{Points stabilized by codimension one subgroups}

Next we compute the connected components of the sets $M^H$, where
$H$ is a codimension one subgroup of $T$. Let $X$ be one of these
components. Then $X^T \neq \emptyset$. Moreover, since $M$ is
simply connected, it is orientable, and hence every connected
component of $M^H$ is orientable. So, if $X$ is not an isolated
point of $M^H$, then it has to be either a circle, a 2-torus, or a
2-sphere, and the first two possibilities are ruled out by the
condition $X^T \neq \emptyset$. We conclude:

\begin{theorem}
Let $H$ be a codimension one subgroup of $T$ and let $X$ be a
connected component of $M^H$. Then $X$ is either a point or a
2-sphere.
\end{theorem}

\begin{remark*} By the Korn-Lichtenstein theorem, every faithful action of
$S^1$ on the 2-sphere is diffeomorphic  to the standard action of
``rotation about the $z-$axis". Therefore the action of the circle
$S^1=T/H$ on the 2-sphere $X$ in the theorem above has to be
diffeomorphic to the standard action. In particular, $\#X^T=2$.
\end{remark*}

We now explicitly determine what these 2-spheres are. By
Theorem~\ref{th:new2}, each of these 2-spheres is the conjugate by
an element of $N_G(T)$ of a 2-sphere containing the identity coset
$p_0 \in M=G/K$; so we begin by determining the 2-spheres
containing $p_0$.

\subsubsection{The space $\fg/\fk$}

The tangent space $T_{p_0}M$ can be identified with $\fg/\fk$, and
the isotropy representation of $T$ on this space decomposes into a
direct sum of two-dimensional $T-$invariant subspaces
\begin{equation}\label{eq:new4}
T_{p_0}M = \oplus V_{[\alpha]} \; ,
\end{equation}
labelled by the roots modulo $\pm 1$,
\begin{equation}\label{eq:new5'}
\alpha \in \Delta_{G,K} / \pm 1 \; .
\end{equation}
One can also regard this as a labelling by the \emph{positive}
roots in $\Delta_{G,K}$; however, since this set of positive roots
is not fixed by the natural action of $W_K$ on $\Delta_{G,K}$,
this is not an \emph{intrinsic} labelling. (This fact is of
importance in Section~\ref{sec:3}, when we discuss the existence
of $G-$invariant almost complex structures on $M$.) Now let $H$ be
a codimension one subgroup of $T$, let $\fh \subset \ft$ be the
Lie algebra of $H$, and let $M^H$ be the set of $H-$fixed points.
Then
$$T_{p_0}M^H = (T_{p_0}M)^H \; .$$
Hence, if $X$ is the connected component of $M^H$ containing
$p_0$, and if $X$ is not an isolated point, then $(T_{p_0}M)^H$
has to be one of the $V_{[\alpha]}$'s in the sum \eqref{eq:new4}.
Hence the adjoint action of $H$ on $\fg/\fk$ has to leave
$V_{[\alpha]}$ pointwise fixed. However, an element $g=\exp{t}$ of
$T$ acts on $V_{[\alpha]}$ by the rotation
\begin{equation}\label{eq:new5}
\chi_{\alpha}(g) = \begin{pmatrix}
  \cos{\alpha(t)} & -\sin{\alpha(t)} \\
  \sin{\alpha(t)} & \cos{\alpha(t)}
\end{pmatrix} \; ,
\end{equation}
so the stabilizer group of $V_{[\alpha]}$ is the group
\begin{equation}\label{eq:new6}
H_{\alpha} = \{ g\in T \; ; \; \chi_{\alpha}(g)=1 \} \; .
\end{equation}

Let $C(H_{\alpha})$ be the centralizer of $H_{\alpha}$ in $G$ and
let $G_{\alpha}$ be the semisimple component of $C(H_{\alpha})$.
Then $G_{\alpha}$ is either $SU(2)$ or $SO(3)$, and since
$G_{\alpha}$ is contained in $C(H_{\alpha})$, $G_{\alpha}p_0$ is
fixed pointwise by the action of $H$. Moreover, since $G_{\alpha}
\nsubseteq K$, the orbit $G_{\alpha}p_0$ can't just consist of the
point $p_0$ itself; hence
\begin{equation}\label{eq:new7}
G_{\alpha}p_0=X \; .
\end{equation}

The Weyl group of $G_{\alpha}$ is contained in the Weyl group of
$G$ and consists of two elements: the identity and a reflection,
$\sigma=\sigma_{\alpha}$, which leaves fixed the hyperplane
$\ker{\alpha} \subset \ft$, and maps $\alpha$ to $-\alpha$.
Therefore, since $\alpha \not\in \Delta_K$, $\sigma_{\alpha}p_0
\neq p_0$, and hence $p_0$ and $\sigma_{\alpha}p_0$ are the two
$T-$fixed points on the 2-sphere \eqref{eq:new7}.

Now let $p=wp_0$ be another fixed point of $T$, with $[w] \in
W_G/W_K$. Let $a$ be a representative for $w$ in $N_G(T)$ and let
$L_a : G \to G$ be the left action of $a$ on $G$. If $X$ is the
2-sphere \eqref{eq:new7}, then the 2-sphere $L_a(X)$ intersects
$M^T$ in the two fixed points $wp_0$ and $w\sigma_{\alpha}p_0$,
and its stabilizer group in $T$ is the group
\begin{equation}\label{eq:new8}
aH_{\alpha}a^{-1}=wH_{\alpha}w^{-1} = H_{w\alpha} \; ,
\end{equation}
where $H_{\alpha}$ is the group \eqref{eq:new6}.

\subsubsection{The GKM graph of $M$}

This concludes our classification of the set of 2-spheres in $M$
which are stabilized by codimension one subgroups of $T$. Now note
that if $X$ is such a two-sphere and $H$ is the subgroup of $T$
stabilizing it, then the orbit space $X/T$ consists of two
$T-$fixed points and a connected one dimensional set of orbits
having the orbitype of $T/H$. Thus these $X$'s are  in one-to-one
correspondence with the edges of the GKM graph of $M$. Denoting
this graph by $\Gamma$ we summarize the graph-theoretical content
of what we've proved so far:

\begin{theorem}\label{th:new2'}
The GKM data associated to the action of $T$ on the homogeneous
space $M=G/K$ is the following.
\begin{itemize}
\item[(1)] The vertices of $\Gamma$ are in one-to-one correspondence
with the elements of $W_G/W_K$;
\item[(2)] Two vertices $[w]$ and $[w']$ are on a common edge of
$\Gamma$ if and only if $[w']=[w\sigma_{\alpha}]$ for some $\alpha
\in \Delta_{G,K}$;
\item[(3)] The edges of $\Gamma$ containing the vertex $[w]$ are in
one-to-one correspondence with the roots, modulo $\pm 1$, in the
set $\Delta_{G,K}$;
\item[(4)] If $\alpha$ is such a root, then the stabilizer group
\eqref{eq:1.1} labelling the edge corresponding to this root is
the group \eqref{eq:new8}.
\end{itemize}
\end{theorem}

In particular, the labelling \eqref{eq:1.1} of the graph $\Gamma$
can be viewed as a labelling by elements $[\alpha]$ of
$\Delta_G/\pm1$. We call this labelling a \emph{pre-axial
function}.

\subsubsection{The connection on $\Gamma$}

One last structural component of the graph $\Gamma$ remains to be
described: Given any graph, $\Gamma$, and vertex, $p$, of
$\Gamma$, let $E_p$ be the set of \emph{oriented} edges of
$\Gamma$ with initial vertex $p$. A \emph{connection} on $\Gamma$
is a function which assigns to each oriented edge, $e$, a
bijective map
$$\theta_{e} : E_p \to E_q \; , $$
where $p$ is the initial vertex of $e$ and $q$ is the terminal
vertex. The graph $\Gamma$ described in Theorem~\ref{th:new2'} has
a natural such connection. Namely, let $e$ be the oriented edge of
$\Gamma$ joining $[w]$ to $[w\sigma_{\alpha}]$. If $e'\in E_{[w]}$
is the oriented edge joining $[w]$ to $[w\sigma_{\beta}]$, then
let $\theta_e(e') = e''$, where $e''$ is the edge joining
$[w\sigma_{\alpha}]$ to $[w\sigma_{\alpha}\sigma_{\beta}]$. This
connection is compatible with the pre-axial function
\eqref{eq:1.1} in the sense that, for every vertex $p$, and every
pair of oriented edges, $e, e' \in E_p$, the roots labelling
$e,e'$, and $e''=\theta_e(e')$ are coplanar in $\ft^*$.

\subsubsection{Simplicity}

A graph is said to be \emph{simple} if every pair of vertices is
joined by at most one edge. Most of the graphs above don't have
this property. There is however an important class of subgroups,
$K$, for which the graph associated with $G/K$ does have this
property.

\begin{theorem}\label{th:simplegraph}
If $K$ is the stabilizer group of an element of $\ft$, then the
graph $\Gamma$ is simple.
\end{theorem}

\begin{proof}
A root $\alpha \in \Delta_G$ is in $\Delta_K$ if and only if the
restriction of $\alpha$ to the subspace $\ft^{W_K}$ of $\ft$ is
zero. Let $\alpha, \beta \in \Delta_{G,K}$ such that $\alpha \neq
\pm \beta$, and let $\sigma_{\alpha}, \sigma_{\beta}$ be the
reflections of $\ft$ defined by $\alpha$ \and $\beta$. Then
$\sigma_{\alpha} \neq \sigma_{\beta}$ and the subspace of $\ft$
fixed by $\sigma_{\alpha}\sigma_{\beta}$ is the codimension 2
subspace on which both $\alpha$ and $\beta$ vanish. If
$\sigma_{\alpha} \sigma_{\beta} \in W_K$, then this subspace
contains $\ft^{W_K}$, so $\alpha$ and $\beta$ are both vanishing
on $\ft^{W_K}$, contradicting our assumption that $\alpha, \beta
\not\in \Delta_K$.
\end{proof}

Another way to prove Theorem \ref{th:simplegraph} is to observe
that $M=G/K$ is a coadjoint orbit of the group $G$. In particular,
it is a Hamiltonian $T-$space and $\Gamma$ is the one-skeleton of
its moment polytope.

\subsection{The GKM definition of the cohomology ring}
\label{ssec:2.3}

We recall how the data encoded in the GKM graph determines the
equivariant cohomology ring $H^*_T(M)$. The inclusion $i : M^T
\to M$ induces a map in cohomology
$$i^* : H^*_T(M) \to H^*_T(M^T) = \text{Maps}(M^T, \SS(\ft^*)) =
\text{Maps}(W_G/W_K , \SS(\ft^*))\; ,$$
and the fact that $M$ is equivariantly formal implies that $i^*$ is
injective. Let $H_T^*(\Gamma)$ be the set of maps
\begin{equation}\label{eq:2.8}
f : W_G/W_K \to \SS(\ft^*)
\end{equation}
that satisfy the compatibility condition:
\
\begin{equation}\label{eq:2.9}
  f([w\sigma_{\alpha}])-f([w]) \in (w\alpha) \SS(\ft^*) \; .
\end{equation}
for every edge $([w],[w\sigma_{\alpha}])$ of $\Gamma$.

The Goresky, Kottwitz and MacPherson theorem \cite{GKM} asserts that
$$
  H^*_T(M) \simeq i^*(H^*_T(M)) = H^*_T(\Gamma) \; .
$$
In the next section we construct a direct isomorphism between this
ring $H^*_T(M)$ and the Borel ring given in \eqref{eq:2.1}.

\subsection{Equivalence between the Borel picture and the GKM picture}
\label{ssec:2.4}

 From the inclusion, $i$, of $M^T$ into $M$, one gets
a restriction map
\begin{equation}\label{eq:2.10}
  i^* : H^*_T(M) \to H^*_T(M^T) \; ;
\end{equation}
and, since $M$ is equivariantly formal, $i^*$ maps $H^*_T(M)$
bijectively onto the subring $H^*_T{\Gamma}$ of $H^*_T(M^T)$.
However, as we pointed out is Section \ref{ssec:2.1},
$$
  H^*_T(M) \simeq \SS(\ft^*)^{W_K}\otimes_{\SS(\ft^*)^{W_G}}
  \SS(\ft^*) \; ;
$$
so, by combining \eqref{eq:2.10} and \eqref{eq:2.1}, we get an
isomorphism
\begin{equation}\label{eq:2.11}
  \K : \SS(\ft^*)^{W_K}\otimes_{\SS(\ft^*)^{W_G}}
  \SS(\ft^*) \to H^*_T({\Gamma}) \; .
\end{equation}

The purpose of this section is to give an explicit formula for
this map. Note that since $M^T$ is a finite set,
$$
  H^*_T(M^T) = \bigoplus_{p \in M^T} H^*_T(p) =  \bigoplus_{p \in M^T}
  \SS (\ft^*) = \text{Maps}(M^T, \SS(\ft^*)) \; .
$$

\begin{theorem}\label{th:4.888.1}
On decomposable elements, $f_1 \otimes f_2$, of the product
\eqref{eq:2.1},
\begin{equation}\label{eq:2.12}
  \K (f_1 \otimes f_2) = g \in \text{Maps}(M^T, \SS(\ft^*)) \; ,
\end{equation}
where, for $w \in W_G$ and $p=wp_0 \in M^T$,
\begin{equation}\label{eq:2.13}
  g(wp_0) = (wf_1) f_2 \; .
\end{equation}
\end{theorem}

\begin{proof}
We first show that \eqref{eq:2.12} and \eqref{eq:2.13} do define a
ring homomorphism of the ring \eqref{eq:2.1} into $H^*(\Gamma)$.
To show that \eqref{eq:2.13} doesn't depend on the representative
$w$ chosen, we note that if $wp_0=w'p_0$, then $\sigma=w(w')^{-1}
\in W_K$. Thus
$$ g(w'p_0)=(w'f_1)f_2 =(w \sigma f_1)f_2 = (wf_1)f_2 = g(wp_0) \;,$$
since $f_1 \in \SS(\ft^*)^{W_K}$. Next, we note that if $f \in
\SS(\ft^*)^{W_G}$, then
$$ \K (f_1f \otimes f_2) = \K (f_1 \otimes f f_2) \; , $$
since
$$w(f_1f)f_2 = (w f_1)(w f) f_2 = (w f_1)ff_2 \; . $$
Thus, by the universality property of tensor products, $\K$ does
extend to a mapping of the ring \eqref{eq:2.1} into the ring
$\text{Maps}(M^T, \SS(\ft^*))$. Next, let $\alpha$ be a root and
let $\sigma \in W_G$ be the reflection that interchanges $\alpha$
and $-\alpha$ and that is the identity on the hyperplane
$$\fh = \{ \xi \in \ft \; ; \; \alpha(\xi)=0 \} \; . $$

Suppose that $p$ and $p'$ are two adjacent vertices of $\Gamma$
with $p'=\sigma p$. To show that $g= \K (f_1 \otimes f_2)$ is in
$H^*_T(\Gamma)$, we must show that the quotient
$$ \frac{g(p')-g(p)}{\alpha} $$
is in $\SS(\ft^*)$. However, if $p=wp_0$, then
$$g(p') - g(p) = (\sigma w f_1 -w f_1)f_2 \; , $$
and since $\sigma$ is the identity on $\fh$, the restriction of
the polynomial $w\sigma f_1$ to $\fh$ is equal to the restriction
of the polynomial $w f_1$ to $\fh$; hence
$$ \frac{g(p')-g(p)}{\alpha}  \in \SS(\ft^*) \; . $$

Finally, we show that the map $\K$ defined by \eqref{eq:2.12} and
\eqref{eq:2.13} has the same equivariance properties with respect
to the action of the Weyl group $W_G$ as does the map
\eqref{eq:2.11}. Note that under the identification
\eqref{eq:2.1}, the action of $W_G$ on $H^*_T(M)$ becomes the
action
$$
  w(f_1 \otimes f_2) = f_1 \otimes w f_2 \; ,
$$
since in the right hand side of \eqref{eq:2.1}, the first factor
is $H^*_G(M)$, so $W_G$ acts trivially on it. In particular,the
ring of $W_G-$invariants in $H^*_T(M)$ is
$$\SS(\ft^*)^{W_K} \otimes_{\SS(\ft^*)^{W_G}} \SS (\ft^*)^{W_G} =
\SS(\ft^*)^{W_K}\; , $$
which is consistent with the identifications
\begin{equation}\label{eq:2.14}
  H^*_G(M)=\SS (\fk^*) = \SS(\ft^*)^{W_K} = H^*_T(M)^{W_G} \; .
\end{equation}

On the other hand, the action of $W_G$ on the space
$$H^*_T(M^T) = \text{Maps}(M^T, \SS(\ft^*)) $$
is just the action
$$
  (wg)(p) = w(g(w^{-1}p)) \; ;
$$
so to check that the map $\K$ defined by \eqref{eq:2.12} and
\eqref{eq:2.13} is $W_G-$equivariant, we must show that if
$$g = \K (f_1 \otimes f_2) \quad \text{ and } g^w = \K (f_1 \otimes wf_2) \; ,$$
then for all points $p=\sigma p_0$,
$$g^w (p) = (wg)(p)\; . $$

However,
$$g^w(p) = (\sigma f_1) (w f_2) = w((w^{-1}\sigma
f_1)f_2)=wg(w^{-1}p) = (wg)(p) \; . $$

Let us now prove that the map $\K$ coincides with the map
\eqref{eq:2.11}. We first note that $\K$ is a morphism of
$\SS(\ft^*)-$modules. For $f \in \SS(\ft^*)$,
$$\K(f_1 \otimes f_2 f) = \K(f_1 \otimes f_2) f\; . $$
Thus, it suffices to verify that $\K$ agrees with the map
\eqref{eq:2.11} on elements of the form $f_1 \otimes 1$.
That is, in view of the identification \eqref{eq:2.14}, it
suffices to show that $\K$, restricted to $\SS(\ft^*)^{W_K}
\otimes 1$, agrees with the map \eqref{eq:2.11}, restricted to
$H^*_T(M)^{W_G}$. However, if $f \in H^*_T(M)^{W_G}$, then $i^* f
\in H^*_T(M^T)^{W_G}$, so it suffices to show that $i^* f$ and
$\K(f \otimes 1)$ coincide at $p_0$, the identity coset of
$M=G/K$. This is equivalent to showing that in the diagram below
$$\xymatrix{
H^*_G(M) \ar[r]\ar[d] & H^*_K(M) \ar[r] &H^*_K(p_0)\ar[d] \\
\SS(\fk^*)^{W_K} \ar[rr] & &  \SS(\fk^*)^{W_K} }$$
the bottom arrow is the identity map. However, the bottom arrow is
clearly the identity on $\SS^0(\fk^*)^K = \CC$ and the two maps on
the top line are $\SS(\fk^*)^K-$module morphisms.
\end{proof}

\section{Almost complex structures and axial functions}
\label{sec:3}

\subsection{Axial functions}
\label{ssec:3.1}

A $G-$invariant almost structure on $M=G/K$ is determined by an
almost complex structure on the tangent space $T_{p_0}M$,
$$J_{p_0} : T_{p_0}M \simeq \fg / \fk \to \fg / \fk \; . $$
For an arbitrary point $gp_0 \in M$, the almost complex structure
on
$$T_{gp_0}M = (dL_g)_{p_0}(T_{p_0}M)$$
is given by
$$J_{gp_0} ((dL_g)_{p_0}(X))=(dL_g)_{p_0}(J_{p_0}(X))\; , $$
for all $X \in \fg / \fk$. This definition is independent on
the representative $g$ chosen if and only if $J_{p_0}$ is
$K-$invariant. Therefore $G-$invariant almost complex structures
on $G/K$ are in one to one correspondence to $K-$invariant almost
complex structures on $\fg / \fk$.

If $M=G/K$ has a $G-$invariant almost complex structure, then the
isotropy representations of $T$ on $T_{p_0}M$ is a complex
representation, and therefore its weights are
well-defined (not just well-defined up to sign). Let
$$T_{p_0}M = \fg/\fk = \bigoplus_{[\beta]} V_{[\beta]}$$
be the root space decomposition of $\fg/\fk$. Then $V_{[\beta]}$
is a one-dimensional complex representation of $T$; let
$\widetilde{\beta} \in \{ \pm \beta \}$ be the weight of this
complex representation:
$$\exp{t} \cdot X_{\widetilde{\beta}} = e^{i \widetilde{\beta}(t)}
X_{\widetilde{\beta}} \quad , \text{ for all } t \in \ft \; . $$

Thus, the map
\begin{equation}\label{eq:3.1}
  s : \Delta_{G,K}/{\pm 1} \to \Delta_{G,K} \; , \;
  s([\beta])=\widetilde{\beta}\; ,
\end{equation}
is a $W_K$-equivariant right inverse of the projection
$\Delta_{G,K}\to \Delta_{G,K}/\{\pm 1\}$. Let $\Delta_0 \subset
\Delta_{G,K}$ be the image of $s$.

The existence of a map \eqref{eq:3.1} is equivalent to the
condition
\begin{equation}\label{eq:3.2}
  w\alpha \neq  - \alpha \quad , \; \text{ for all } w \in W_K \; , \;
  \alpha \in \Delta_{G,K} = \Delta_G - \Delta_K \; ,
\end{equation}
hence \eqref{eq:3.2} is a necessary condition for the existence of
a $G-$invariant almost complex structure on $M$. We will see in
the next section that this condition is also sufficient.

We can now define a labelling of the \emph{oriented} edges,
$E_{\Gamma}$, of the GKM graph $\Gamma$, as follows. Let $[w]\in
W_G/W_K$ be a vertex of the graph and let
$e=([w],[w\sigma_{\beta}])$ be an oriented edge of the graph, with
$\beta \in \Delta_0$. This edge corresponds to the subspace
$V_{[w\beta]}$ (see \eqref{eq:new8}) in the decomposition
$$T_{[w]}M = \bigoplus_{\beta \in \Delta_0} V_{[w\beta]} \; , $$
and the $G-$invariance of the almost complex structure implies
that $T$ acts on $V_{[w\beta]}$ with weight $w\beta$. We define
$\alpha: E_{\Gamma} \to \ft^*$ by
\begin{equation}\label{eq:3.3}
\alpha([w],[w\sigma_{\beta}]) = w\beta \; , \text{ for all } \beta
\in \Delta_0, w \in W_G \; .
\end{equation}

\begin{theorem}\label{th:axialfcn}
The map $\alpha : E_{\Gamma} \to \ft^*$ has the following
properties:
\begin{enumerate}
\item If $e_1$ and $e_2$ are two oriented edges with the same
initial vertex, then $\alpha(e_1)$ and $\alpha(e_2)$ are linearly
independent;

\item If $e$ is an oriented edge and $\bar{e}$ is the same edge,
with the opposite orientation, then $\alpha(\bar{e}) = -
\alpha(e)$;

\item If $e$ and $e'$ are oriented edge with the same initial
vertex, and if $e''=\theta_e(e')$, then $\alpha(e'')-\alpha(e')$
is a multiple of $\alpha(e)$.
\end{enumerate}
\end{theorem}

\begin{proof}
The first assertion is a consequence of the fact that the only
multiples of a root $\alpha$ that are roots are $\pm \alpha$.

If $e$ is the oriented edge that joins $[w]$ to $[w\sigma_{\beta}]$
and that is labelled by $w\beta \in w\Delta_0$, then
$$\alpha(\bar{e}) = (w\sigma_{\beta})(\beta)= -w\beta =
-\alpha(\bar{e}) \; . $$

Finally, if $e$ joins $[w]$ to $[w\sigma_{\beta}]$ and if $e'$
joins $[w]$ to $[w\sigma_{\gamma}]$ (with $\beta, \gamma \in
\Delta_0$), then $e''$ joins $[w\sigma_{\beta}]$ to
$[w\sigma_{\beta}\sigma_{\gamma}]$, and
$$\alpha(e'')-\alpha(e)= w\sigma_{\beta}\gamma -w\gamma =
w(\sigma_{\beta}\gamma-\gamma) = -\langle \gamma, \beta \rangle
w\beta = - \langle \gamma, \beta \rangle \alpha(e) \; . $$
\end{proof}

Equivalently, Theorem~\ref{th:axialfcn} says that
$\alpha: E_{\Gamma} \to \ft^*$ is an \emph{axial function}
compatible with the connection $\theta$, in the sense of \cite{GZ1}.

\subsection{Invariant almost complex structures}
\label{ssec:3.2}

As we have seen in Section~\ref{ssec:3.1}, \eqref{eq:3.2} is a
necessary condition for the existence of a $G-$invariant almost
complex structure on $M=G/K$; in this section we show that it is
also a sufficient condition.

\begin{theorem}\label{th:converse}
If the condition
$$w\alpha \neq  - \alpha \quad , \; \text{ for all } w \in W_K \; ,
\;
  \alpha \in \Delta_{G,K} = \Delta_G - \Delta_K \; ,
$$
is satisfied, then $M$ admits a $G-$invariant
almost complex structure.
\end{theorem}

\begin{proof}
Consider the complex representation of $K$ on $(\fg/\fk)_{\CC} =
\fg_{\CC} / \fk_{\CC}$ and let
$$(\fg/\fk)_{\CC} = \bigoplus_{j} V_j$$
be the decomposition into irreducible representations;
$(\fg/\fk)_{\CC}$ is self dual, hence
$$\bigoplus_{j} V_j =(\fg/\fk)_{\CC} = (\fg/\fk)_{\CC}^* = \bigoplus_{j}
V_j^*=\bigoplus_{j} \overline{V_j} \; $$
Therefore $\overline{V_j} = V_l$ for some $l$. If $\alpha$ is a
highest weight of $V_j$, then condition \eqref{eq:3.2} implies
that $-\alpha$ is not a weight of $V_j$; however, $-\alpha$ is a
weight of $\overline{V_j}$, hence $\overline{V_j} \neq V_j$.
Therefore
$$(\fg/\fk)_{\CC} = \bigoplus_{j} (V_j \oplus \overline{V_j}) = U
\oplus \overline{U} $$
as complex $K-$representations, and this induces a $K-$invariant
almost complex structure
$$J : \fg / \fk \to \fg / \fk$$
as follows: If $x \in \fg / \fk$, then there exists a unique $y
\in \fg / \fk$ such that $x+iy \in U$, and we define $J(x)=y$. As
we have shown before, this is equivalent to the existence of a
$G-$invariant almost complex structure on $M$.
\end{proof}

An alternative way of proving Theorem~\ref{th:converse} is to
observe that the condition \eqref{eq:3.2} is equivalent to the
existence of a $W_K-$equivariant section $s: \Delta_{G,K}/{\pm 1}
\to \Delta_{G,K}$. Let $s$ be such a section and let $\Delta_0
\subset \Delta_G - \Delta_K$ be the image of $s$. Then (see
\eqref{eq:new4})
$$\fg / \fk = \bigoplus_{\alpha \in \Delta_0} V_{[\alpha]}$$
and one can define a $K-$invariant almost complex structure $J$ by
requiring that for each $\alpha \in \Delta_0$, $J$ acts on
$V_{[\alpha]}$ by
\begin{equation}\label{eq:3.4}
J\begin{pmatrix}
  X_{\alpha} \\
 X_{-\alpha}
\end{pmatrix} = \begin{pmatrix}
  X_{-\alpha} \\
 -X_{\alpha}
\end{pmatrix}\; .
\end{equation}

\section{Morse theory on the GKM graph}
\label{sec:4}

\subsection{Betti numbers}
\label{ssec:4.1}

Henceforth we assume that $M$ admits a $G-$invariant almost
complex structure, determined (see \eqref{eq:3.4}) by the image
$\Delta_0 \subset \Delta_{G,K}$ of a section $s : \Delta_{G,K} /
{\pm 1} \to \Delta_{G,K}$. Let $\Gamma$ be the GKM graph of $M$
and let
$$\alpha : E_{\Gamma} \to \ft^*$$
be the axial function \eqref{eq:3.3}. Then the edges whose initial
vertex is the identity coset in $W_G/W_K$ are labelled by vectors in
$\Delta_0$.

Let $\xi \in \ft$ be a regular element of $\ft$, \emph{i.e.}
$$\beta(\xi) \neq 0 \quad , \quad \text{ for all } \beta \in
\Delta_G \subset \ft^* \; . $$
For a vertex $[w] \in W_G/W_K$, define the \emph{index} of $[w]$
to be
$$
ind_{[w]} = \# \{ e \in E_{[w]} \; ; \; \alpha(e)(\xi ) <0 \} \; ,
$$
and for each $k \geq 0$, let the $k-$th \emph{Betti number} of
$\Gamma$ be defined by
$$
  \beta_{k}(\Gamma) = \# \{ [w] \in  W_G/W_K \; ; \; ind_{[w]}=k
  \} \; .
$$
The index of a vertex obviously depends on $\xi$, but the Betti
numbers do not.

\begin{theorem}\cite{GZ1}
The Betti numbers $\beta_{k}(\Gamma)$ are combinatorial invariants
of $\Gamma$ (\emph{i.e.} are independent of $\xi$).
\end{theorem}

In general these Betti numbers are not equal to the Betti numbers
$$\beta_{2k}(M) = \dim H^{2k}(M)$$
of $M=G/K$; however, we show in the next
section that there is a large class of homogeneous spaces for
which they are equal. One should note
that $\beta_{2k}(M)$ is the dimension of the \emph{ordinary}
cohomology of $M$ as a \emph{vector space}, while
$\beta_k(\Gamma)$ counts the number of generators of degree $2k$ of
the \emph{equivariant} cohomology ring of $M$, as a \emph{free module}
over the symmetric algebra $\SS(\ft^*)$.

\subsection{Morse functions}
\label{ssec:4.2}

Let $\xi \in \ft$ be a regular element.

\begin{definition}
A function $f : W_G/W_K \to \RR$ is called a \emph{Morse function}
compatible with $\xi$ if for every oriented edge $e=([w],[w'])$ of
the GKM graph, the condition $f([w']) > f([w])$ is satisfied whenever
$\alpha(e)(\xi)>0$.
\end{definition}

Morse function do not always exist; however, there is a simple
necessary and sufficient condition for the existence of a Morse
function: Every regular element $\xi
\in \ft$ determines an orientation $o_{\xi}$ of the edges of
$\Gamma$: an edge $e \in E_{\Gamma}$ points upward
(with respect to $\xi$) if $\alpha_e(\xi) >0$, and points
downward if $\alpha_e(\xi)<0$. The associated directed graph $(\Gamma,
o_{\xi})$ is the graph with all upward-pointing edges.

\begin{proposition}
There exists a Morse function compatible with $\xi$ if and only if
the directed graph $(\Gamma, o_{\xi})$ has no cycles.
\end{proposition}

\subsection{Invariant complex structures}
\label{ssec:4.3}

In this section we show that the existence of Morse functions on
the GKM graph, which is a combinatorial condition, has geometric
implications for the space $M=G/K$.

\begin{theorem}
The GKM graph $(\Gamma, \alpha)$ admits a Morse function
compatible with a regular $\xi \in \ft$ if and only if the
almost complex structure determined by $\alpha$ is a $K-$invariant
complex structure on $M$.  Moreover, if this is the case, then the
combinatorial Betti numbers agree with the topological Betti numbers.
That is,
$$b_{k}(\Gamma) = b_{2k}(M) \quad . $$
\end{theorem}

\begin{proof}
Let $f : W_G/W_K \to \RR$ be a Morse function compatible with
$\xi$, and let $[w]$ be a vertex of the GKM graph where $f$
attains its minimum. If we replace $\xi$ by $w^{-1}(\xi)$ and $f$
by $(w^{-1})^*f$, then the minimum of this new function is $p_0$.
Thus, without loss of generality, we may assume that the minimum
vertex $[w]$ is the identity coset in $W_G/W_K$. Then
$$\Delta_0 = \{ \beta \in \Delta_{G,K} \; ; \; \beta(\xi) > 0 \} \; , $$
hence $\Delta_0$ is the intersection of $\Delta_{G,K}$ with the
positive Weyl chamber determined by $\xi$. Let
$$\fp = \fk_{\CC} \oplus (\bigoplus_{\beta \in \Delta_0}
\fg_{\beta}) \; . $$
Then $\fp$ is a parabolic subalgebra of $\fg_{\CC}$, hence the
almost complex structure determined by $\alpha$ is actually a
complex structure.

If $G_{\CC}$ is the simply connected Lie group with Lie algebra
$\fg_{\CC}$ and if $P$ is the Lie subgroup of $G_{\CC}$
corresponding to $\fp$, then
$$M = G / K = G_{\CC} / P \; ,$$
hence $M$ is a coadjoint orbit of $G$. Then $M$ is a Hamiltonian
$T-$space and the GKM graph of $M$ is the 1-skeleton of the moment
polytope, and therefore the combinatorial Betti numbers agree with
the topological Betti numbers.

On the other hand, if the almost complex structure is integrable
then $\fp$ is a parabolic subalgebra of $\fg_{\CC}$, $M=G/K
\subset \fg^*$ is a coadjoint orbit of $G$, and for a generic
direction $\xi \in \ft \subset \fg$, the map $f : W_G/W_K \to \RR$
given by
$$f([w]) = \langle [w], \xi \rangle$$
(with $W_G /W_K \to G/K \to \fg^*$) is a Morse function on the GKM
graph compatible with $\xi$.
\end{proof}

\section{Examples}
\label{sec:5}

\subsection{Non-existence of almost complex structures}
\label{sec:5.1}

Let $G$ be a compact Lie group such that $\fg_{\CC}$ is the simple
Lie algebra of type $B_2$.
\begin{figure}[h]
\begin{center}
\includegraphics{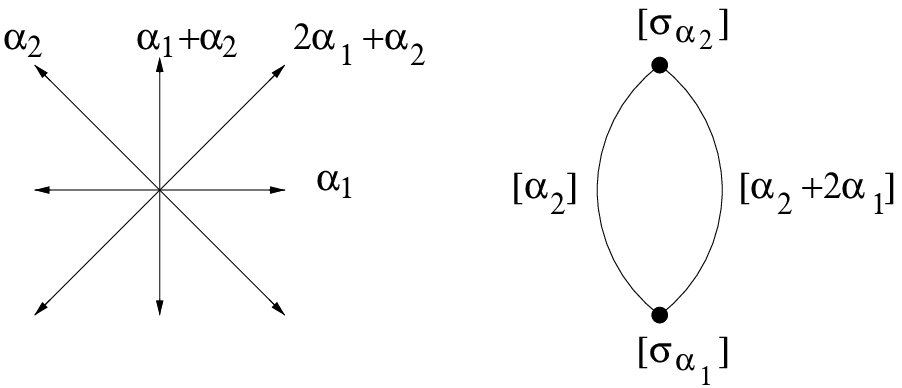}
\caption{} \label{fig:b2-d2}
\end{center}
\end{figure}
Let $\alpha_1, \alpha_1+\alpha_2$ be the short positive roots and
let $\alpha_2, \alpha_2+2\alpha_1$ be the long positive roots. Let
$K$ be the subgroup of $G$ corresponding to the root system
consisting of the short roots. Then $\fk_{\CC} = D_2 = A_1 \times
A_1$ and $K \simeq SU(2) \times SU(2)$. The quotient $W_G/W_K$ has
two classes: the class of $\sigma_{\alpha_1} \in W_K$ and the
class of $\sigma_{\alpha_2} \in W_G - W_K$.

The GKM graph $\Gamma$ has two vertices, joined by two edges, and
the edges are labelled by $[\alpha_2] , [\alpha_2+2\alpha_1] \in
\Delta_{G,K}/{\pm 1}$. If
$w=\sigma_{\alpha_1+\alpha_2}\sigma_{\alpha_1} \in W_K$, then
$w\alpha_2=-\alpha_2$ and $\alpha_2 \in \Delta_{G,K}$, hence one
can't define an axial function on $\Gamma$.  In this example,
$G/K=S^4$, which is does not admit an almost complex structure.

\subsection{Non-existence of Morse functions}
\label{sec:5.2}

Let $G$ be a compact Lie group such that $\fg_{\CC}$ is the simple
Lie algebra of type $G_2$. Let $\alpha_1, \alpha_1+\alpha_2$, and
$2\alpha_1 + \alpha_2$ be the short positive roots and let
$\alpha_2, 2\alpha_2+3\alpha_1, \alpha_2+3\alpha_1$ be the long
positive roots. Let $K$ be the subgroup of $G$ corresponding to
the root system consisting of the short roots. Then $\fk_{\CC} =
A_2$ and $K \simeq SU(3)$. The quotient $W_G/W_K$ has two classes:
the class of $\sigma_{\alpha_1} \in W_K$ and the class of
$\sigma_{\alpha_2} \in W_G - W_K$.
\begin{figure}[h]
\begin{center}
\includegraphics{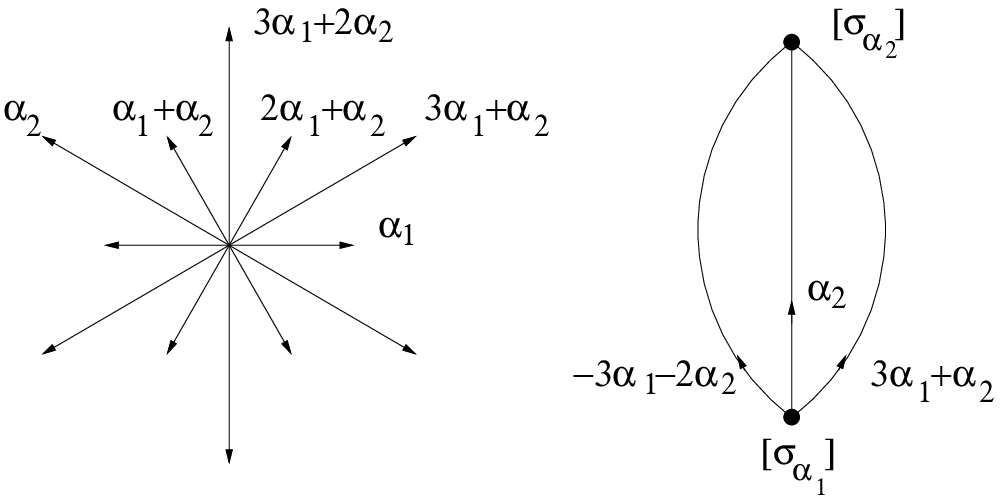}
\caption{} \label{fig:g2-a2}
\end{center}
\end{figure}

The GKM graph $\Gamma$ has two vertices, joined by three edges,
and the edges are labelled by $[\alpha_2] , [2\alpha_2+3\alpha_1],
[\alpha_2+3\alpha_1] \in \Delta_{G,K}/{\pm 1}$. There are two
$W_K-$equivariant sections of the projection $\Delta_{G,K} \to
\Delta_{G,K}/{\pm 1}$, corresponding to $\{ \alpha_2,
\alpha_2+3\alpha_1,-2\alpha_2-3\alpha_1\}$ and $\{ -\alpha_2,
-\alpha_2-3\alpha_1,2\alpha_2+3\alpha_1\}$. If
$$\Delta_0=\{ \alpha_2,\alpha_2+3\alpha_1,-2\alpha_2-3\alpha_1\}\; , $$
then the axial function is shown in Figure~\ref{fig:g2-a2} and
there is no Morse function on $\Gamma$: the corresponding almost
complex structure is not integrable.  In this example, $G/K=S^6$,
which admits an almost complex structure, but no invariant complex
structure.

\subsection{The existence of several almost complex structures}
\label{sec:5.3}

Let $G=SU(3)$ and $K=T$. Then the homogeneous space $G/K$ is the
manifold of complete flags in $\CC^3$. The root system of $G$ is
$A_2$, with positive roots $\alpha_1, \alpha_2$, and
$\alpha_1+\alpha_2$ of equal length. The Weyl group of $G$ is
$W_G=S_3$, the group of permutations of $\{1,2,3\}$, and $W_K=1$,
hence $W_G/W_K =W_G = S_3$.

The GKM graph is the bi-partite graph $K_{3,3}$ : it has 6
vertices and each vertex has 3 edges incident to it, labelled by
$[\alpha_1], [\alpha_2]$, and $[\alpha_1+\alpha_2]$. There are
$2^3$ possible $W_K-$invariant sections, hence eight $G-$invariant
almost complex structures on $G/K$. If
$$\Delta_0 = \{\alpha_1, \alpha_2, \alpha_1+\alpha_2 \}\; ,  $$
then the corresponding almost complex structure is integrable and
there is a Morse function on $\Gamma$ compatible with $\xi \in
\ft$ such that both $\alpha_1(\xi)$, and $\alpha_2(\xi)$ are
positive. This Morse function is given by $f(w) = \ell(w)$ where
$\ell(w)$ is the length of $w$, \emph{i.e.}, in this case, the
number of inversions in $w$.
\begin{figure}[h]
\begin{center}
\includegraphics{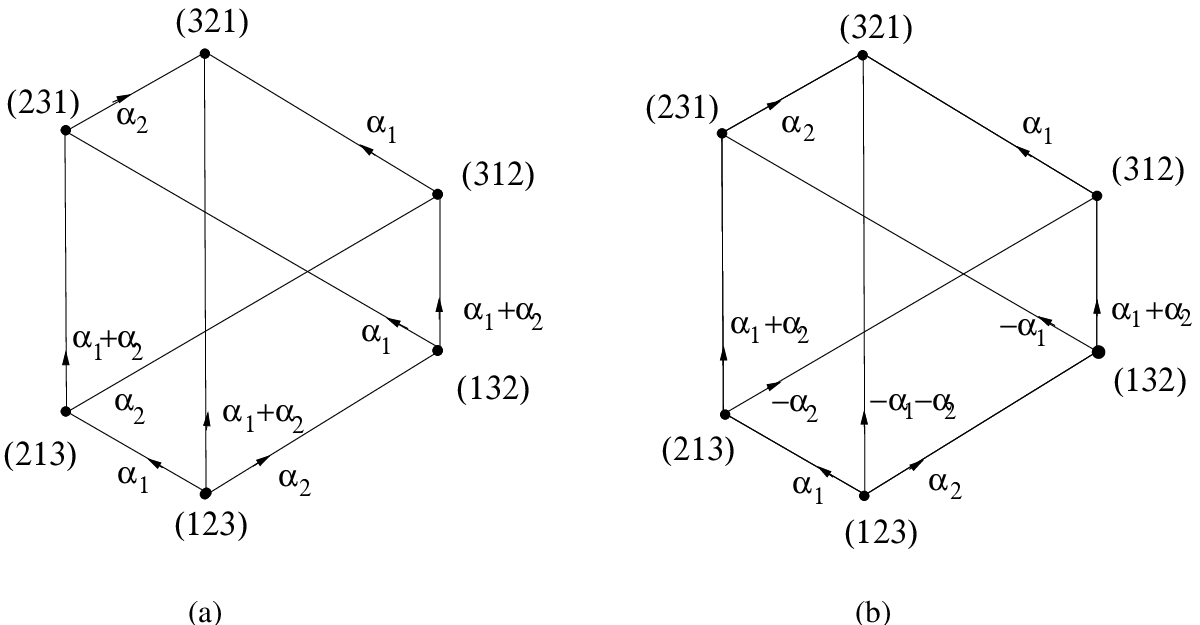}
\caption{} \label{fig:a2}
\end{center}
\end{figure}

However, if
$$\Delta_0 = \{\alpha_1, \alpha_2, -\alpha_1-\alpha_2 \}\; ,  $$
then the corresponding almost complex structure is not integrable
and there is no Morse function on $(\Gamma, \alpha)$ : for every
vertex $w$ of $\Gamma$, there exist three edges $e_1, e_2$, and
$e_3$, going out of $w$, such that
$$\alpha_{e_1}+\alpha_{e_2}+\alpha_{e_3} = 0 \;  ,$$
hence there is no vertex of $\Gamma$ on which a Morse function
compatible with some $\xi \in \ft$ can achieve its minimum.

\end{document}